\documentclass{amsart}

\usepackage{amsfonts,amsthm,amsmath,amssymb,mathrsfs}
\usepackage[hypertex, linktocpage]{hyperref}
\usepackage[all]{xy}

\DeclareMathOperator{\st }{\ \vert \ }

\DeclareMathOperator{\supp}{supp}
\DeclareMathOperator{\Exp}{Exp}
\DeclareMathOperator{\Log}{Log}

\newcommand{\n}{\par\noindent}

\newcommand{\sn}{\par\smallskip\noindent}
\newcommand{\mn}{\par\medskip\noindent}
\newcommand{\bn}{\par\bigskip\noindent}
\newcommand{\EL}{{EL}}
\newcommand{\LE}{{LE}}
\newcommand{\lra}{\longrightarrow}
\newcommand{\eps}{\varepsilon}
\newcommand {\notion }[1]{{\bf #1}}%
\newcommand {\kk }{{\bf k}}%

\def\N{\mathbb{N}}
\def\Z{\mathbb{Z}}
\def\R{\mathbb{R}}

\newtheorem{Thm}{Theorem}[section]
\newtheorem{Cor}[Thm]{Corollary}
\newtheorem{Lem}[Thm]{Lemma}
\newtheorem{Prop}[Thm]{Proposition}
\theoremstyle{definition}
\newtheorem{Ex}[Thm]{Example}
\newtheorem{Rem}[Thm]{Remark}
\newtheorem{Def}[Thm]{Definition}

\title[EL and LE series]{Comparison of Exponential-Logarithmic and Logarithmic-Exponential series}
\author{Salma Kuhlmann}
\address{Fachbereich Mathematik und Statistik, Universit\"at Konstanz,
Germany}
 \email{salma.kuhlmann@uni-konstanz.de}
\author{Marcus Tressl}
\address{University of Manchester, School of Mathematics, UK }
\email{marcus.tressl@manchester.ac.uk}
\date{\today}
\thanks{We would like to thank Micka\"el Matusinski for reading a
preliminary version of this paper and providing helpful comments.}
\begin{document}
\begin{abstract}   We
explain how the field of logarithmic-exponential series constructed
in \cite{DMM1} and \cite {DMM2} embeds as an exponential field in
any field of exponential-logarithmic series constructed in
\cite{KK1}, \cite {K} and \cite {KS}. On the other hand, we explain
why no field of exponential-logarithmic series embeds in the field
of logarithmic-exponential series. This clarifies why the two
constructions are intrinsically different, in the sense that they
produce non-isomorphic models of Th$(\R_{\mbox{an, exp}})$; the
elementary theory of the ordered field of real numbers, with the
exponential function and restricted analytic functions.
\end{abstract}
\maketitle
\tableofcontents
\section{Introduction}
\label{introduction} Several constructions of non-archimedean
ordered fields endowed with an exponential and a logarithmic map,
and with differential operators have appeared in the literature in
the last two decades. All constructions are based on the use of
fields of generalized series (see Definition \ref{defnK((G))}). Let
us summarize some of the literature highlights. The early works of
\cite{D} and \cite{DG} were motivated by Tarski's open problem
concerning the decidability of $T_{\exp}$:= the elementary theory of
$(\R, +, \cdot, <, \exp )$. In \cite{W}, Wilkie established the
model completeness of this theory. It became particularly
interesting to understand the algebraic structure of the
non-archimedean models of $T_{\exp}$. In Dahn, \cite[p.183]{D}, a
(non-surjective) exponential map is defined on a field of
generalized series. This field, denoted by $L_0$, is reconsidered in
the works \cite{DMM}, \cite{DMM1}, \cite{DMM2}. There, the authors
construct a non-archimedean model of $T_{\exp}$, the so-called field
of LE-series. The construction is based on Dahn's field $L_0$, and
is attributed to Dahn-G\"oring \cite{DG} in \cite[p.63]{DMM2}. The
non-standard models constructed in \cite{KK1}, \cite{KK2},
\cite{KKS}, \cite{K}, \cite{KS} are based on so-called {\it
prelogarithmic fields} (see Section 3 ). In \cite{E}, \cite{H1},
\cite{H2}, \cite{KM}, and \cite{KM2} the focus is on the
differential operators and on power series expansions for solutions
of differential equations. Example \ref{multiplicativehahn} and
variants thereof (e.g. with $\Gamma =\{\log_m x\>;\> m\in \N _0\}\>$
instead) are considered in \cite[2.3.1]{S}.

All constructions described in the literature are complicated and
hard to read for the non-specialists. The inherent complexity of the
construction makes a comparison of the resulting structures a tour
de force. In this paper, we endeavor to explain these different
constructions in a unified way.  More precisely, we isolate a list
of main basic construction steps. We explain the constructions of
\cite{DMM1} and \cite {DMM2} and those of \cite{KK1}, \cite {K} and
\cite {KS} in light of this list. We show how one may combine them
in different ways to recover different constructions. This unified
approach has several advantages. On the one hand, it makes the
constructions more accessible, on the other, it allows to compare
the various constructions, up to isomorphism of the obtained models.
\mn The paper is organized as follows: In section
\ref{preliminaries}, we recall some preliminaries on fields of
generalized power series. In sections \ref{prelogs}, \ref{example},
\ref{expext} and \ref{expclosure} we explain the construction of
exponential-logarithmic power series fields $\R((G))^{\mbox{ EL}}$
(starting with an arbitrary ordered abelian group $G$; cf. \cite
{K}). The construction is split into 3 different steps: the initial
data is a given prelogarithmic section (see Section \ref{prelogs}
for the definition), which is used to construct in the second step
the exponential extension, and finally in the third step the
exponential closure. In Section \ref{morphisms}, we introduce the
notion of morphisms in the category of fields endowed with
prelogarithmic sections. Morphisms are used in Section
\ref{section:LE} to ``run'' a generalized construction of
logarithmic-exponential fields in exponential-logarithmic fields. To
get as a special case an embedding of the LE-series field (cf.
\cite{DMM1}, \cite{DMM2}) in any EL-series field (cf. \cite{K}) we
need to make a special choice of $G$. The assumptions on $G$ are
explained and analyzed in Section \ref{section:groups_of_exponents}.
In Section \ref{LEinEL} , we then show how to embed the LE-series
field of \cite{DMM1}, \cite{DMM2} into the exponential closure of
the subfield generated by logarithmic words (cf.\ \cite{H1},
\cite{H2}, \cite{S}), which is in turn a subfield of any EL-series
field. In particular, we get an embedding of \'Ecalle's ``trig\`ebre
de transseries'', which is identified as a subfield of the LE-series
field in \cite [Section 7.1]{DMM2}. The last section \ref{notELinLE}
is devoted to proving that no non-archimedean exponential field $\R
((G))^{\rm EL}$ does embed as an exponential field
 into the LE-series field.

\section{Reminder on generalized power series.} \label{preliminaries}
\noindent Throughout, ${\bf k}$ will denote a totally ordered field
and $(G\>, \cdot\>, 1\>, <)$ a multiplicative totally ordered
abelian group.  \begin{Def}\label{defnK((G))} The field of {\it
generalized power series} over $\> G \>$ with coefficients in ${\bf
k}$ is defined as follows:
\[{\bf k}((G)):=\{ \alpha =\sum \alpha(g)\> g\st \alpha (g) \in \mathbf{k} \>,\> \supp
\alpha \text { is anti well ordered}\} .\] Here $\supp \alpha :=\{
g\in G\st \alpha (g)\neq 0\} $ is the support of $\alpha $.
Addition and multiplication of series are defined in the usual
manner. See \cite{N} for more details. \sn
  The field ${\bf k}((G))$ is endowed with a {\it canonical
valuation } $v: {\bf k}((G)) \rightarrow G\cup \{0\}$, defined by
$v(0)= 0$ where $0 <G$, $v(\alpha):= \max \supp \alpha$ for
$\alpha\not= 0$. The valuation verifies that $v(\alpha \alpha ')\> =
\> v(\alpha)\cdot v(\alpha ')\>,$ as well as the ultrametric
triangle inequality: $v(\alpha + \alpha ') \leq \max \{v(\alpha),
v(\alpha')\}$, with equality if $v(\alpha)\not= v(\alpha')\>.$ Thus
$({\bf k}((G)), v)$ is a {\it valued field} with {\it value group}
$\> G \>$ ($G$ is also called the group of exponents, or the group
of monomials), and {\it residue field} $\>{\bf k}$. The field ${\bf
k}((G))$ comes equipped with the {\it anti lexicographic order}:
$0\not= \alpha \in {\bf k}((G))$ is $>0$ if $\alpha (v(\alpha))
>0$. The valuation $v$ is {\it compatible } with this order, i.e.
for all $\alpha_1, \alpha_2$ with $0 <\alpha_1 <\alpha_2$ we have
$v(\alpha_1) \leq v(\alpha_2)$.
\end{Def}
\sn
 We identify ${\bf k}$ with the ordered subfield
 $\{ \alpha\cdot 1\st
\alpha\in {\bf k}\}\>$ of $\> {\bf k}((G))$. Similarly, we
identify $G$ with the ordered subgroup $\{ 1\cdot g\st g\in G\} $
of ${\bf k}((G))^{\times}\>$. Hence $G>0$ in ${\bf k}((G))$ and
one might think of
 $G$ as a multiplicative group of germs of functions
 ${\bf k} \lra {\bf k} $ at $+\infty $.
 \sn
 For any $g\in G$ we denote by $G^{<g}$ the
strict initial segment of $G$ determined by $g$, similar notations
for $G^{\leq g}$, $G^{>g}$, and $G^{\geq g}$. More generally, for
any subset $A\subseteq G$ we denote by $G^{<A}:=\{g \st g<a,
\mbox{ for all } a \in A \>\}$, similar notations for $G^{\leq
A}$, $G^{>A}$, $G^{\geq A}\>$. For any subset $S \subseteq G$, let
${\bf k}((S)):=\{ \alpha \in {\bf k}((G))\st \supp \alpha
\subseteq S\} $. Since $\supp (\alpha_1+\alpha_2) \subseteq \supp
\alpha_1 \cup \supp \alpha_2$ and $\supp (-\alpha)=\supp \alpha$,
$({\bf k}((S)) , +)$ is a ${\bf k}$-linear subspace of $({\bf
k}((G)) , +)$.\sn
 The {\it valuation ring} ${\bf k}((G^{\leq 1}))={\bf k}\oplus {\bf k}((G^{< 1}))$ is
a convex subring of ${\bf k}((G))$ with maximal ideal ${\bf
k}((G^{<1}))$. If $\alpha \in {\bf k}((G))$, $\alpha >0$ and
$g:=\max \supp \alpha $, then $g^{-1}\cdot \alpha \in {\bf
k}((G^{\leq 1})) $, hence there are uniquely determined $a\in {\bf
k}^{>0}$, $\eps \in {\bf k}((G^{<1}))$ such that $\alpha =g\cdot
a\cdot (1+\eps)$. We have the following direct sum (respectively,
multiplicative direct sum) decompositions:
\begin{equation} \label{decomp1}
{\bf k}((G))= {\bf k}((G^{>1}))\oplus {\bf k}\oplus {\bf k}((G^{<
1}))\>,
\end{equation}
\begin{equation} \label{decomp2}
{\bf k}((G))^{>0}=G\cdot {\bf k}^{>0}\cdot (1+{\bf k}((G^{<1})))
\>.
\end{equation}
\begin{Rem}
In the literature, and in particular in the papers of the first
author cited here, the group $G$ is written additively, $g\in G$ is
associated to the monomial $t^{g}$ and series in ${\bf k}((G))$ are
written $\sum \alpha _g t^g$ and have well ordered support. The map
$g\mapsto t^g$ is an order {\it reversing} isomorphism, the order on
${\bf k}((G))$ is the lexicographic order, and the canonical
valuation is defined by $v(\alpha):=\min \supp \alpha$. In the
present paper, we have opted for the multiplicative notation because
it is more suitable for functional interpretations (cf.\ \cite{K2}
and \cite{KM}).
\end{Rem}
\sn
\section{Prelogarithmic
sections and associated prelogarithms.} \label{prelogs}
 \sn In this and the
following three sections, we shall recall the construction, and
give examples of exponential-logarithmic series fields. We begin
by recalling the definition of the restricted logarithm defined on
the goup of 1-units of the valuation ring of ${\bf k}((G))$. \mn
\begin{Def}
 The {\it logarithm on
1-units} is the map:
\begin{equation} \label{restrictedlog}
 1+{\bf k}((G^{<1}))\lra {\bf k}((G^{< 1})) \mbox{ defined by }
1+\eps \mapsto \sum _{i\geq 1}(-1)^{(i-1)} \> {\eps ^i\over i} \>.
\end{equation}
\end{Def}
\sn This map is an isomorphism of ordered groups (cf. \cite {A}).
\mn We now focus on extending the domain of the logarithm to ${\bf
k}((G))^{>0}$.
\begin{Def}
 A {\it prelogarithmic section} $l$ of
${\bf k}((G))$ is an embedding of ordered groups $l:(G\>,
\cdot)\lra ({\bf k}((G^{>1})),+)\>.$ The tuple $({\bf k}((G)),l)$
is called a {\it series field with prelogarithmic section} $l$;
if $l$ is clear from the context we just say
$k((G))$ is a \textit{prelogarithmic series field}.
\end{Def}
\sn Note that prelogarithmic sections on ${\bf k}((G))$ exist,
cf.\ Remark \ref{Hahnembedding} and Example
\ref{multiplicativehahn} below.
\begin{Def}
Let ${\bf k}$ be an ordered field and $\log :({\bf k}^{>0},\cdot
)\lra ({\bf k},+)$ an order preserving embedding of groups. We say
that $({\bf k}\>, \log)$ is an {\it ordered prelogarithmic field}.
If $\log$ is surjective we say that $({\bf k}\>, \log)$ is an
 {\it ordered logarithmic field}, or that $({\bf k}\>, \exp)$ is an
{\it ordered exponential field}, where $\exp:=\log ^{-1}$.
\end{Def}
\sn Let $({\bf k}\>, \log)$ be an ordered prelogarithmic field and
$G$ be a multiplicative, abelian ordered group. Let $l$ be a
prelogarithmic section of ${\bf k}((G))$. If $\alpha \in {\bf
k}((G))$, $\alpha >0$ and $g:=\max \supp \alpha $, write $\alpha
=g\cdot a\cdot (1+\eps)$ with $a\in {\bf k}^{>0}$, $\eps \in {\bf
k}((G^{<1}))$.
\begin{Def}
Define
\begin{equation}\label{(L)}
L(\alpha)=L(g\cdot a\cdot (1+\eps))=l(g)+\log a +\sum _{i\geq
1}(-1)^{(i-1)} \> {\eps ^i\over i}
\end{equation}
Then $L:({\bf k}((G))^{>0},\cdot )\lra ({\bf k}((G)),+)$ is the
uniquely determined order preserving embedding of groups,
extending $\log $, $l$ and the logarithm on 1-units. We call $L$
the {\it prelogarithm} associated to the prelogarithmic section
$l$ on ${\bf k}((G))$.
\end{Def}
\sn
\section  {Main examples of prelogarithmic series
fields}\label{example} \sn
\begin{Def}
Let $\Gamma$ be a totally ordered set, and ${\bf k}$ a totally
ordered field. We shall begin by defining $ \Gamma ^{{\bf k}}$;
the {\it Hahn group} of rank $\Gamma$ with exponents in ${\bf
k}\>$, written multiplicatively. Consider a totally ordered set of
variables
\[X:= \{\ x_{\gamma}\>\st\> \gamma \in \Gamma\}\] and let $\Gamma^{\bf k}$ be the set
of formal products of the form \[g = \prod _{\gamma \in \Gamma}
x_{\>\gamma} ^{\>g(\gamma)}\] where $g(\gamma)\in {\bf k}$ and
$\supp g:=\{ \gamma \in \Gamma \st g(\gamma)\neq 0\} $ is anti
wellordered in $\Gamma$. Multiplication is defined pointwise:
\[g_1g_2 = \prod _{\gamma \in \Gamma} x_{\>\gamma}^{\>g_1(\gamma)+
g_2(\gamma)}\>\] so $1$ is the product with empty support. The
order is anti lexicographic: $g
> 1$ if $g(\gamma) > 0$ where $\gamma = \max \supp g$. For example
$x_{\gamma} > 1$ for all $\gamma\in \Gamma$. Thus, $\Gamma^ {\bf
k}$ is a totally ordered abelian group.
\end{Def}
\begin{Rem} \label{Hahnembedding}
Hahn's Embedding Theorem \cite{H} states that every totally
ordered abelian group G is (isomorphic to) a subgroup of a Hahn
group $\Gamma^{ \R}$, for a suitably chosen $\Gamma$. Therefore
when constructing prelogarithmic series fields $\R((G))$ one may assume
without loss of generality that $G$ is a subgroup of a Hahn group.
\end{Rem}
\begin{Ex} \label{multiplicativehahn}
 Set $G=\Gamma^
{\bf k}$, we define the {\it basic prelogarithmic section} on
${\bf k}((G))$ by :
\[l(\prod _{\gamma \in \Gamma} x_{\>\gamma} ^{\>g(\gamma)}) = \sum
_{\gamma \in \Gamma} g(\gamma) x_{\gamma}\>.\] \mn Now assume as
above that ${\bf k}$ admits a  $\log :({\bf k}^{>0},\cdot )\lra
({\bf k},+)$, and let $L$ be the associated {\it basic
prelogarithm} as given by (\ref{(L)}). We denote by ${\bf
k}((\Gamma^{\bf k}))^L$ the prelogarithmic series field thus
constructed.
\end{Ex}
\mn
 Our aim is to use
prelogarithmic series fields to construct ordered exponential
fields $(K, E)$ which are models of $T_{\exp}$. So we are mainly
interested in exponentials satisfying the {\it growth axiom}
scheme: \mn {\bf (GA)}\ \ $\alpha \in K$, $\alpha \geq
n^2\;\Longrightarrow\; E(\alpha)> \alpha ^n \hspace{4cm} (n\geq
1)$ \mn Note that because of the hypothesis $\alpha \geq n^2$,
{\bf (GA)} is only relevant for $v(\alpha) \geq 1$. Let us
consider the case $v(\alpha)>1$. In this case,``$\alpha>n^2\,$''
holds for all $n\in\N$ if $\alpha$ is positive, and the axiom
scheme {\bf (GA)} is thus equivalent to the assertion
\begin{equation}                            \label{f(x)>x^n}
\forall n\in\N:\; E(\alpha)>\alpha ^n \>.
\end{equation}
Applying $L:= E^{-1}$ on both sides, we find that this is
equivalent to
\begin{equation}                            \label{increll}
\forall n\in\N:\; \alpha > L(\alpha ^n)=nL(\alpha) \>.
\end{equation}
Via the valuation $v$, this in turn holds if
\begin{equation}                            \label{stronglog}
v(L(\alpha)) < v(\alpha) \>.
\end{equation}
\begin{Def}
A (pre)logarithm $L$ on the ordered field $K$ will be called a
{\bf (GA)}-(pre)logarithm if it satisfies
(\ref{stronglog}) for $\alpha \in K^{>0}$ with $v(\alpha)>1$.
\end{Def}
\begin{Ex} \label{sigmamultiplicativehahn}
  The basic prelogarithm $L$ on ${\bf k}((G))$ (given in Example \ref {multiplicativehahn})
  does {\it not} satisfy
  {\bf (GA)} (e.\ g. $L(x_{\gamma})=x_{\gamma}$ ). To
remedy to this problem, we fix an embedding \[\sigma\>: \Gamma
\longrightarrow \Gamma\] which is decreasing (i.e. $\sigma(\gamma)
< \gamma$ for all $\gamma \in \Gamma$), and define the induced
prelogarithmic section $l_{\sigma}$ as follows:
\[l_{\sigma}(\prod _{\gamma \in \Gamma} x_{\>\gamma} ^{\>g(\gamma)}) = \sum
_{\gamma \in \Gamma} g(\gamma) x_{\sigma(\gamma)}\>.\] The
associated prelogarithm ( given in (\ref{(L)}) ) is denoted by
$L_{\sigma}$. We note that $L_{\sigma}$ satisfies {\bf (GA)}: one
verifies that (\ref{stronglog}) holds for $\alpha \in {\bf
k}((G))^{>0}$, with $g:= v(\alpha) >1$. Set $x_{\gamma} :=\max
\supp g\>\>;$ one verifies that $$x_{\sigma(\gamma)} =
v(l_{\sigma} (g)) = v(L_{\sigma}(\alpha)) < v(\alpha)=g \mbox{ in
}G^{>1}\>.$$ \sn We denote by ${\bf k}((\Gamma^{\bf k}))^{\sigma
L}$ the prelogarithmic series field thus constructed.
\end{Ex}
\begin{Ex} \label{Schmeling}
In Example \ref{multiplicativehahn} above, $\Gamma$ can consist of
a totally ordered set of germs at infinity of non-oscillating real
valued functions of a real variable, which are infinitely large
and positive (i.e. $\mbox{ lim }_{x\rightarrow +\infty} f(x) =
+\infty$). For example take
\[\Gamma:=\{\log_m x\>\st\> m\in \mathbb Z\}\>.\]
Here, $\log$ is the natural logarithm on $\R$ and $\log _m$ is its
m-th iterate (for $n \in \N\>,$ $\log _{-n}= \exp_n$ the n-th
iterate of the exponential). Let
\[\sigma: \log_m x \mapsto \log_{m+1} x \>.\]
Let $G$ be the group of {\it exponential-logarithmic words}; that
is, the subgroup of $\Gamma^{ \R}$ consisting of products with
{\it finite} support. As in Example \ref{multiplicativehahn}
above, the induced prelogarithmic section $l_{\sigma}$ on
$\R((G))$ is:
\[l_{\sigma}(x^{r_0}(\log x)^{r_1}\cdots (\log_n x)^{r_n}) = r_0 \log x +
\cdots + r_n \log _{n+1} x.\]
\end{Ex}
\sn
\section{The exponential extension of a prelogarithmic
field.}\label{expext}
 \sn
 The
following construction or close variants of it, namely to add
exponentials to a given field with prelogarithmic section, is used
by \cite {D}, \cite {E}, \cite {DMM1}, \cite {DMM2}, \cite {H1},
\cite {H2}, \cite {KK1}, \cite {KK2}, \cite{K}, \cite {KS} and \cite
{S}. It is the core step in constructing exponentials of infinitely
large elements. \sn From now on, we fix a prelogarithmic series
field $({\bf k}((G)),l)$. Observe that $l:G\lra ({\bf
k}((G^{>1})),+)\>$ is not surjective if $G\not=\{1\}$ (cf.\
\cite{KKS}). More precisely, the exponential of any element in ${\bf
k}((G^{>1}))\setminus l(G)$ is not defined. We shall enlarge our
group of monomials $G$ to a group extension $G^\#$ to include the
missing exponentials.\sn We take $G^\#$ to be a {\it multiplicative}
copy $e[{\bf k}((G^{>1}))]$ of ${\bf k}((G^{>1}))$ over $l(G)$. More
precisely, we construct
 $G^\#$ formally as follows:
 $$G^\#:=\{e(\alpha)\>\>\st\>\>
\alpha \in {\bf k}((G^{>1})), \mbox{ where } e(\alpha):=g \mbox{
if } \exists g\in G \mbox{ s.t. } \alpha =l(g) \>\}$$ By its
definition, $G$ is a subset of $G^\#$.
 We define multiplication on $G^\#$ as follows:
  $e(\alpha_1)e(\alpha_2):=e(\alpha_1 + \alpha_2 )$.
 In particular, if $g_1=e(\alpha_1)$, $g_2=e(\alpha_2) \in G\>,$ then $e(\alpha_1)e(\alpha_2)=
 e(l(g_1) + l(g_2)) = e(l(g_1g_2))=g_1g_2\>,$ so $G$ is a {\it subgroup} of $G^\#\>$.
 We equip $G^\#$ with a total order: $e(\alpha_1) < e(\alpha_2)$
 if and only if $\alpha_1 < \alpha_2$ in ${\bf k}((G^{>1}))\>.$ Again, if $g_1=e(\alpha_1)$,
 $g_2=e(\alpha_2) \in G\>,$ then $e(\alpha_1) < e(\alpha_2)$ if
 and only if $l(g_1) < l(g_2)$ in ${\bf k}((G^{>1}))$ if and only if $g_1 < g_2$ in $G$,
  so $G$ is an {\it ordered
subgroup} of $G^\#\>.$  Since $G\subseteq G^\#$ as ordered abelian
multiplicative groups, we view ${\bf k}((G))$ as an ordered
subfield of ${\bf k}((G^\#))$ by identifying ${\bf k}((G))$ with
the elements of ${\bf k}((G^\#))$ having support in $G$. \sn One
verifies that the map $l^\#:(G^\#,\cdot )\lra ({\bf k}((G^{\#\
>1})),+)$ defined by $l^\#(e(\alpha)):=
\alpha$ for $\alpha \in {\bf k}((G^{>1}))\>$ is a prelogarithmic
section with $l^\#(G^\#) = {\bf k}((G^{>1}))$ and $l^\#$ extends
$l$ on $G$. By construction of the logarithms $L$ and $L^\#$ on
${\bf k}((G))^{>0}$ and ${\bf k}((G^\#))^{>0}$ respectively,
$L^\#$ is an extension of $L$.
\begin{Def} We define the {\it exponential extension} of
$({\bf k}((G)),L)$ to be $({\bf k}((G^\#)),L^\#)$.
\end{Def}\mn We have the following commutative diagram
of embeddings (whenever we have a set theoretic
inclusion we use the arrow $\hookrightarrow $):
\begin{equation}\label{(figure:expextension)}
\vcenter{
\xymatrix@=40pt{
\mathbf{k}((G^\#))^{>0}\ar [rr]^{L^\#}&&\mathbf{k}((G^\#))\\
G^\#\ar@{^{(}->}[u] \ar [r]^{l^\#}\ar [rd]_{\cong }^{l^\#}&\mathbf{k}((G^{\#>1}))\ar@{^{(}->}[ru]&\\
G\ar [r]_l\ar@{^{(}->} [u]&\mathbf{k}((G^{>1}))\ar@{^{(}->} [u]\ar@{^{(}->}[r]&\mathbf{k}((G))\ar@{^{(}->}[uu]
}}
\end{equation}
\sn
\section{The exponential closure of a prelogarithmic field}
\label{expclosure}
 \sn
From now on we assume that the given logarithm on the residue
field $\log :({\bf k}^{>0},\cdot )\lra ({\bf k},+)$ is surjective.
\begin{Def} If $n=0$ set $({\bf k}((G))^{\#
n},L^{\# n}):=({\bf k}((G)),L)\>.$ For $n\in \N\>,$ define
inductively the {\it $n$-th exponential extension} of $({\bf
k}((G)),L)$: $({\bf k}((G))^{\# n},L^{\# n}):={\rm  the\
exponential\ extension\ of\ }({\bf k}((G^{\# n-1})),L^{\#
n-1})\>.$
\end{Def} \sn Hence ${\bf k}((G))^{\# n}={\bf k}((G^{\# n}))$.
Here, the prelogarithm on ${\bf k}((G^{\# n}))^{>0}$ associated to
$l^{\# n}$ is denoted by $L^{\# n}$.
\begin{Def} Set ${\bf {\bf k}((G))^{\EL}}:=\bigcup _{n\in \N _0}\>
{\bf k}((G))^{\# n}$ and $ \Log:=\bigcup _{n\in \N _0}\> L^{\#
n}\>$. We call $({\bf {\bf k}((G))^{\EL}}\>, \Log)$ is the
EL-series field over the prelogarithmic field $({\bf k}((G))\,
l)$.
\end{Def} \sn Below, we gather some properties of $({\bf
k}((G))^{\EL}\>,\Log)$.
 We note that ${\bf
k}((G))$ is contained in the image of $L^\#$: the image of $L^\#$
contains ${\bf k}((G^{<1}))$ (from the logarithm on 1-units),
${\bf k}((G^{>1}))$ (as image of $l^\#$) and ${\bf k}$, since
$\log $ is surjective. By induction we have
\begin{equation} \label{loggeneral1}
{\bf k}((G^{\# n}))\subseteq L^{\# n+1}[{\bf k}((G^{\#
n+1}))^{>0}]
\end{equation}
and
\begin{equation} \label{loggeneral2}
{\bf k}((G^{\# n\>,\> >1})) = L^{\# n+1}(G^{\# n+1})\>.
\end{equation}
 So
$\Log:({\bf k}((G))^{\EL\>,\>
>0},\cdot )\lra ({\bf k}((G))^\EL ,+)$ is an order preserving
{\it isomorphism}. Let ${\bf \Exp} :({\bf k}((G))^\EL ,+)\lra
({\bf k}((G))^{\EL \>,\> >0},\cdot )$ denote the inverse of
$\Log$. By (\ref {loggeneral1}) and (\ref {loggeneral2}), we have
\begin{equation}\label{exp general}
\Exp[{\bf k}((G^{\# n}))]\subseteq {\bf k}((G^{\# n+1}))^{>0}
\end{equation}
and
\begin{equation}
 \Exp[{\bf k}(((G^{\#
n\>,\> >1}))]=G^{\# n+1}.
\end{equation}
Finally, although we do not use this fact, we note that ${\bf
k}((G))^\EL \subseteq {\bf k}((G^\EL ))$ where $G^\EL :=\bigcup
_nG^{\# n}$.
\begin{Ex} \label{multiplicativehahn2}
We consider the prelogarithmic field ${\bf k}((\Gamma^{\bf k}))^L$
constructed in Example \ref{multiplicativehahn}. The EL series
field obtained as above by forming its exponential closure shall
be henceforth denoted by ${\bf k}((\Gamma^{\bf k}))^{\EL}$.
\end{Ex}
\begin{Ex} \label{sigmamultiplicativehahn2}
We consider the prelogarithmic field ${\bf k}((\Gamma^{\bf
k}))^{\sigma L}$ constructed in Example
\ref{sigmamultiplicativehahn}. The EL series field obtained as
above by forming its exponential closure shall be henceforth
denoted by ${\bf k}((\Gamma^{\bf k}))^{\sigma\EL }$. We claim that
$\Log _{\sigma}\>:\>({\bf k}((\Gamma^{\bf k}))^{\sigma\EL\>,\>
>0},\cdot )\lra ({\bf k}((\Gamma^{\bf k}))^{\sigma\EL} ,+)$ satisfies {\bf
(GA)}. In Example \ref{sigmamultiplicativehahn}, we have shown
that\begin{equation} \label{ga}
 v(l_{\sigma}(g)) < g \mbox{ for } g\in
G^{>1}\>. \end{equation} It suffices to show that this property is
preserved by exponential extension. Let $g^\# \in G^{\#\>,\> >1}$.
Set $g:=v(l_{\sigma}^{\#}(g^\#))\in G^{>1}\>.$ Appying (\ref{ga})
we get:\[ v(l_{\sigma}(v( l_{\sigma}^{\#}(g^\#)))) <
v(l_{\sigma}^{\#}(g^\#))\>.\] Since $l_{\sigma}(v(
l_{\sigma}^{\#}(g^\#))) > 0$ and $ l_{\sigma}^{\#}(g^\#) >0\>,$ we
get that (by compatibility of $v$ with the order):
\[l_{\sigma}(v( l_{\sigma}^{\#}(g^\#)))<
l_{\sigma}^{\#}(g^\#)\>.\] Since $l_{\sigma}^\#$ extends
$l_{\sigma}$ , the last inequality reads
\[l_{\sigma}^\#(v( l_{\sigma}^{\#}(g^\#)))<
l_{\sigma}^{\#}(g^\#)\>.\] Since $l_{\sigma}^\#$ is order
preserving, we conclude  that $v( l_{\sigma}^{\#}(g^\#)))<
g^\#\>,$ as required.
\end{Ex}
\section {Morphisms of prelogarithmic fields.} \label{morphisms}
\sn In this section, we define a morphism  on a prelogarithmic
series field $({\bf k}((G)), l)$ induced by a morphism $\psi$ of
$G$, and use it in turn to induce a special automorphism (`` value
group induced'' automorphism) $\psi^{\EL}$ of the corresponding
EL-series field $({\bf k}((G))^{\EL}, \Log)$. These automorphisms
$\psi^{\EL}$ will play the role of ``$\Log$-substitution maps''
(cf.\ Example \ref{sigmamultiplicativehahn3}); symbolically
$\psi^{\EL}(f(x))=f(\Log (x))$ (substituting $\Log (x)$ for $x$).
The aim is to recover in Section \ref {LEinEL} the inverse map
$(\psi^{\EL})^{ -1}$ (substituting $\Exp (x)$ for $x$) which plays
a key role in the construction of the LE series field ${\bf
k}((x^{-1}))^{\LE}$ of \cite{DMM2}.
\begin{Def}
Let $\nu : {\bf k}((G_1))\lra {\bf k}((G_2))$ be a homomorphism of
${\bf k}$--vector spaces. We say that $\nu$ {\it respects
arbitrary sums} if it commutes with arbitrary sums, i. e. if $\nu
(\sum a_gg)=\sum a_g\nu (g)\>$.
\end{Def}
\sn Let $({\bf k}((G_1)),l_1)$ and $({\bf k}((G_2)),l_2)$ be
series fields with prelogarithmic sections and let $\psi :G_1\lra
G_2$ be an order preserving embedding (of groups). Then $\psi $
extends to ${\bf k}((G_1))$ in the obvious way and the resulting
extension is also denoted by $ \psi $. Since $\psi$ respects
arbitrary sums (from the way it is defined), $ \psi$ commutes with
the logarithms on 1-units defined in (\ref{restrictedlog}).
\begin{Def} We call $\psi $ a {\it morphism} from $({\bf
k}((G_1)),l_1)$ to $({\bf k}((G_2)),l_2)$ if $\psi $ commutes with
the prelogarithmic sections, that is: $\> \psi \circ l_1=l_2\circ
\psi \>$ on $G_1$:
\begin{equation}\label{(figure:morphism)}
\vcenter{
\xymatrix@=40pt{
\mathbf{k}((G_1))\ar [r]^{ \psi }&\mathbf{k}((G_2))\\
\mathbf{k}((G_1^{>1}))\ar [r]\ar@{^{(}->}[u]&\mathbf{k}((G_2^{>1}))\ar@{^{(}->}[u]\\
G_1\ar [r]_\psi \ar[u]^{l_1}&G_2\ar[u]_{l_2}}}
\end{equation}
If this is the case, then $\psi $ also commutes with the
prelogarithms: $$\psi \circ L_1=L_2\circ \psi \> \mbox{ on }\>
{\bf k}((G_1))^{>0}\>.$$ (Indeed, if $a\in {\bf k}((G_1))^{>0}$
write $a=g\cdot r\cdot (1+\eps )$ and compute: \n $\psi
(L_1(a))=\psi (l_1g+\log r+\log (1+\eps ))= \psi l_1(g)+\log
r+\psi (\log (1+\eps ))=l_2\psi(g)+\log r+\log (1+\psi \eps )=
L_2(\psi (a))\>$). \mn If $\psi $ is such a morphism, then $\psi$
induces a morphism
$$\psi ^\#\>:\> ({\bf k}((G_1^\# )),l_1^\# )\longrightarrow ({\bf
k}((G_2^\# )),l_2^\# )$$ defined as follows: for $g^\# \in G_1^\#
\>$ $$\psi ^\# (g^\# )=(l_2^\#)^{-1}\circ \psi\circ l_1^\# (g^\#
)\>.$$ Observe that $\psi\circ l_1^\# (g^\# )\in {\bf
k}((G_2^{>1}))$ is indeed in the image of $l_2^\# $, and that
$\psi ^\#$ extends $\psi$.
\end{Def}
\sn Here is the situation in a diagram. The part of the diagram
with non broken arrows, commutes: \bn
\begin{equation}\label{(figure:morphismsmall)}
\vcenter{
\xymatrix@=20pt{
\mathbf{k}((G_1^\# ))\ar[rrrrr]^{\psi ^\# }&&&&&\mathbf{k}((G_2^\# ))&&\\
&&G_1^\#\ar[rrrrr]^<(0.25){\psi ^\# =(l_2^\#)^{-1}\circ \psi \circ
l_1^\#}
\ar@{_{(}->}[ull]\ar@{-->}[dl]^{l_1^\#}_\cong &&&&&G_2^\#\ar@{_{(}->}[ull]\ar@{-->}[dl]^{l_2^\#}_\cong\\
&\mathbf{k}((G_1^{>1}))\ar@{^{(}->}[dl]&&&&&\mathbf{k}((G_2^{>1}))\ar@{^{(}->}[dl]&\\
\mathbf{k}((G_1))\ar'[rr][rrrrr]^{ \psi} \ar@{^{(}->}[uuu]&&&&&\mathbf{k}((G_2))\ar@{^{(}->}'[uu][uuu]&&\\
&&G_1\ar[rrrrr]^\psi \ar@{^{(}->}[ull]\ar@{^{(}->}[uuu]&&&&&G_2\ar@{_{(}->}[ull]\ar@{^{(}->}[uuu]
}}
\end{equation}
\bn Adding the prelogarithmic sections of diagram \ref
{(figure:expextension)} to diagram \ref {(figure:morphismsmall)}
we get a complete illustration of the constructions so far in the
following figure: \bn
\begin{equation}\label{(figure:morphismbig)}
\vcenter{
\xymatrix@=13pt{
\mathbf{k}((G_1^\# ))\ar[rrrrr]^{\psi ^\# }&&&&&\mathbf{k}((G_2^\# ))&&\\
&&\mathbf{k}((G_1^\#))^{>0}\ar@{-->}[ull]_<(0.15){L_1^\#}&&&&&
\mathbf{k}((G_2^\#))^{>0}\ar@{-->}[ull]_<(0.15){L_2^\#}\\
&&G_1^\#\ar@{^{(}->}[u]\ar[rrrrr]^<(0.25){\psi ^\#
=(l_2^\#)^{-1}\circ \psi\circ
l_1^\#}\ar@{^{(}->}[uull]\ar@{-->}[dl]^{l_1^\#}_\cong
&&&&&G_2^\#\ar@{^{(}->}[u]\ar@{_{(}->}[uull]\ar@{-->}[dl]^{l_2^\#}_\cong\\
&\mathbf{k}((G_1^{>1}))\ar@{^{(}->}[dl]&&&&&\mathbf{k}((G_2^{>1}))\ar@{^{(}->}[dl]&\\
\mathbf{k}((G_1))\ar'[rr][rrrrr]^{ \psi} \ar@{^{(}->}[uuuu]&&&&&\mathbf{k}((G_2))
\ar@{^{(}->}'[uu][uuuu]&&\\
&&&&&&\\
&&G_1\ar@{-->}[uuul]_<(0.25){l_1}\ar[rrrrr]^\psi \ar@{^{(}->}[uull]\ar@{_{(}->}[uuuu]
&&&&&G_2\ar@{-->}[uuul]_<(0.25){l_2}\ar@{_{(}->}[uull]\ar@{_{(}->}[uuuu]
}}
\end{equation}
 We now iterate the $\sharp $-construction: by induction we get a
morphism $$\psi ^{\# n}\>:\> ({\bf k}((G_1^{\# n})),l_1^{\#
n})\longrightarrow ({\bf k}((G_2^{\# n})),l_2^{\# n})$$ extending
$\psi ^{\# n-1}$.  In the following two lemmas, we record for
later use properties of the resulting map
$$\psi ^{\EL }\>:\> {\bf k}((G_1))^{\EL } \lra {\bf k}((G_2))^{\EL } \>,$$
defined by $$\psi ^{\EL }:=\bigcup_{n\in \N} \psi^{\# n}\>.$$
\begin{Prop}\label{properties of psi}
 The following properties hold.
 \sn
 (1) \
$\psi ^{\EL }$ is a ${\bf k}$--embedding of ordered fields which
respects arbitrary sums .\sn (2) \ $\psi ^{\EL } $ respects $\Log$
and $\Exp$: for $\> a\in {\bf k}((G_1))^{\EL }$; $\psi ^\EL (\Log
(a))=\Log (\psi ^\EL (a))\>$ (if $a>0$) and
  $\psi ^{\EL } (\Exp (a))=\Exp (\psi ^\EL (a))\>.$  \sn (3) \ If $G_2\subseteq \psi ^\EL (G_1^{\# n})$
for some $n\in \N $ then ${\bf k}((G_2^{\# m}))$ is in the image
of ${\bf k}((G_1^{\# n+m}))$ under $\psi ^\EL $ for all $m\in \N
$. Hence in this case $\psi ^\EL $ is surjective and $(\psi ^\EL
)^{-1}$ respects arbitrary sums.
\end{Prop}
\begin{proof}
(1) and (2) are inherited from the corresponding properties of
$\psi ^\# $. To see (3), note that by definition, ${\bf k}((G_2))$
is in the image of $\psi ^\EL \vert _{{\bf k}((G_1^{\# n}))}$.
Thus ${\bf k}((G_2^{\# m}))$ is in the image of $\psi ^\EL \vert
_{{\bf k}((G_1^{\# n+m}))}$ for all $m\in \N $.
\end{proof}
\begin{Def}
Assume that $G=G_1=G_2$ and that $l=l_1=l_2$. We say that the
morphism $\psi: G \lra G$ is a {\it contracting} morphism if for
all $g\in G^{>1}\>$ $\psi (g)< g\>.$ Similarly, we say that $\psi
^\EL$ is a  {\it contraction} if $\psi ^\EL (g)<g$ for all $g\in
G^\EL $, $g>1$.
\end{Def}
\sn Contracting morphisms will play a key role in
section \ref{LEinEL}. We record the following for later use.

\begin{Prop}\label{psi contracting} Assume that $G=G_1=G_2$ and that
$l=l_1=l_2$, and that $\psi$ is a contracting morphism, then $\psi
^\EL$ is a contraction
\end{Prop}
\begin{proof}
By induction, it is enough to show $\psi ^\# (g^\# )< g^\# $ for
all $g^\# \in G^\# $, $g^\# >1$. Since $\psi (g)<g$ ($g\in G$,
$g>1$) we have $ \psi (f)<f$ for all  $f\in {\bf k}((G^{>1}))$,
$f>0$.  Take $g^\# \in G^\#, g^\# >1$. Then $l^\#(g^\#)\in {\bf
k}((G^{>1}))$ is positive, Thus $ \psi (l^\#(g^\#))<l^\#(g^\#)$
and since $(l^\#)^{-1}$ is monotone, we get $\psi ^\# (g^\#
)=(l^\#)^{-1}\circ \psi\circ l^\#(g^\#)<
(l^\#)^{-1}(l^\#(g^\#))=g^\# $ as desired.
\end{proof}
\begin{Ex} \label{sigmamultiplicativehahn3}
We study the morphism induced by $\sigma$ on the
exponential-logarithmic series field ${\bf k}((\Gamma^{{\bf
k}}))^{\sigma\EL}$ examined in Examples
\ref{sigmamultiplicativehahn} and \ref{sigmamultiplicativehahn2}.
The embedding $\sigma:\Gamma \lra \Gamma$ lifts to an order
preserving embedding $\psi_{\sigma}:=\tilde{\sigma}: G \lra G$ of
$G:=\Gamma^ {\bf k}$ into itself, defined by $$\psi_{\sigma}(\prod
_{\gamma \in \Gamma} x_{\>\gamma} ^{\>g(\gamma)}) \>=\>\prod
_{\gamma \in \Gamma} x_{\>\sigma(\gamma)} ^{\>g(\>\gamma)}=\>\prod
_{\gamma \in \Gamma} l_{\>\sigma}(x_{\gamma}) ^{\>g(\>\gamma)}.$$
It is straightforward to verify that $\psi_{\sigma}$ is a
morphism, which extends to $\psi^{\>\sigma \EL}\>:\>{\bf
k}((\Gamma^{\bf k}))^{\sigma\EL}\lra {\bf k}((\Gamma^{\bf
k}))^{\sigma\EL}\>.$ Note also that for all $g\in G^{>1}\>$
$\psi_{\sigma} (g)< g\>,$ so that $\psi ^{\>\sigma\EL} (g)<g$ for
all $g\in G^\EL $, $g>1$ (by Proposition \ref{psi contracting}).
We also apply Proposition \ref{properties of psi} (3) to establish
that $\psi ^{\>\sigma\EL}$ is surjective (even if $\sigma$ is
not): it suffices to show that for all $g\in G$ we have that $g\in
\psi_{\sigma}^\#(G^\#)$. Let $g=\prod _{\gamma \in \Gamma}
x_{\>\gamma} ^{\>g(\gamma)} \in G$. Now $\sum_{\gamma \in \Gamma}
g(\gamma)\>x_{\gamma} \in {\bf k}((G^{>1})) =
l_{\sigma}^\#(G^\#)$. Let $g^\# \in G^\#$ be such that
$l_{\sigma}^\#(g^\#)= \sum_{\gamma \in \Gamma}
g(\gamma)\>x_{\gamma}\>.$ We compute: $\psi_{\sigma}^\#
(g^{\#})\>=\> l_{\sigma}^{\#\> -1} \circ \psi_{\sigma} \circ
l_{\sigma}^\# (g^\#)\>=\> \> l_{\sigma}^{\#\> -1} \circ
\psi_{\sigma} (\sum_{\gamma \in \Gamma}
g(\gamma)\>x_{\gamma})\>=\>l_{\sigma}^{\#\> -1}(\sum_{\gamma \in
\Gamma} g(\gamma)\>x_{\sigma(\gamma}))\>=\> g$ as required.
\end{Ex}
\section {Groups satisfying the growth axiom.}
\label{section:groups_of_exponents}
\begin{Lem}\label{characterize lex product}
Let $A,B$ be subgroups of the ordered abelian group $(C,+)$ and
suppose $A<\vert b\vert $ for all $b\in B$, $b\neq 0$. Then $A$ is
a convex subgroup of  $A\oplus B$.
\end{Lem}
\begin{proof}
Clearly $A+B$ is a direct sum. Let $a,a'\in A$ and $b\in B$ with
$0<a+b<a'$. From $0<a+b$ we get $-b<a$. From $a+b<a'$ we get
$b<a'-a$. Hence $A\not <\vert b\vert $ and our assumption implies
$b=0$.\end{proof} \sn We return to ${\bf k}((G))$ with
prelogarithmic section $l:G\lra {\bf k}((G^{>1}))$, prelogarithm
$L:({\bf k}((G))^{>0},\cdot )\lra ({\bf k}((G)),+)$, and
associated EL-series field ${\bf k}((G))^\EL $ with logarithm
$\Log $ and exponential $\Exp $.

\medskip
\begin{Lem}\label{characterize local growth}
Let $U\subseteq H\subseteq G$ be subgroups such that $U$ is a
proper convex subgroup of $H$ ($U=\{1\} $ is not excluded).
Suppose  $$\Log(H)<\vert f\vert \text { for all }f\in {\bf
k}((H^{>U})),\ f\neq 0.\leqno{(\dagger)}$$ Let
$$H^{\# ,U}=H\cdot \Exp({\bf k}((H^{>U})))\subseteq G^{\#}.$$
Then $H$ is a proper convex subgroup of $H^{\# ,U}$, $H^{\#
,U}$ is the anti-lexicographic product of $H$ and $\Exp({\bf
k}((H^{>U})))$ and
 $$\Log(H^{\# ,U})<\vert f\vert \text { for all }f\in {\bf k}(((H^{\# ,U})^{>H})),\ f\neq 0.$$
\end{Lem}
\begin{proof}
$H$ is a proper subgroup of $H^{\# ,U}$, since for every $h\in
H$ with $h>U$ we know $H<\Exp (h)$ from $(\dagger )$. Let $h^\#
\in H^{\# ,U}$, hence $h^\# =h\cdot \Exp (f_1)$ with $h\in
H$ and $f_1\in {\bf k}((H^{>U}))$. By assumption there is $h_1\in
H$ with $U<h_1$. Then $\Log (h^\# )=\Log (h)+f_1<h_1+f_1\in
{\bf k}((H))$, hence for every $f\in {\bf k}(((H^{\#
,U})^{>H}))$ with $f>0$ we have $\Log (h^\# )\leq h_1+f_1<f$.
Finally Lemma \ref {characterize lex product} tells us that
condition $(\dagger )$ implies that $H^{\# ,U}$ is the
anti-lexicographic product of $H$ and $\Exp({\bf
k}((H^{>U})))$.\end{proof}

\begin{Def}\label{defn:growth_axiom}
Let $H$ be a subgroup of $G$ with $H\neq \{ 1\} $.
We say that \notion{the growth axiom holds for }$H$ if condition
$(\dagger )$ of Lemma \ref {characterize local growth} is satisfied for
$H$ and $U=\{ 1\} $, in other words if
$$\Log(H)<\vert f\vert \text { for all }f\in {\bf k}((H^{>1})),\ f\neq 0.$$
\end{Def}
Observe that for a simply generated subgroup $H=h^\Z $ with $h>1$,
$H$ satisfies the growth axiom if and only if
$$v(\Log (h))<v(h)$$
where $v$ denotes the natural valuation (in multiplicative
notation), i.e. the valuation whose valuation ring is the convex
hull of $\Z $. This is what classically is called ``the growth
axiom for $h$".

\smallskip
We iterate \ref {characterize local growth} to obtain

%
\begin{Prop}\label{growth axiom for H}
Let $({\bf k},\log)$ be a logarithmic field and let $({\bf
k}((G)),l)$ be a field with prelogarithmic section.
 Let $H\neq \{ 1\} $ be a subgroup of $G$, satisfying the growth axiom
(cf. Definition \ref {defn:growth_axiom}). Let $H^n$ be defined
inductively by $H^{-1}=\{ 1\} $, $H^0=H$ and
$$H^{n+1}=H^n\cdot \Exp {\bf k}(([H^n]^{>H^{n-1}}))\subseteq G^{\# n+1}.$$
Then $H^n$ is a proper convex subgroup of $H^{n+1}$,
$H^{n+1}$ is the anti-lexico\-graphical product of
$H^n$ and $\Exp {\bf k}(([H^n]^{>H^{n-1}}))$
and
$$\Log (H^{n+1})<|f| \text { for all }f\in \kk (((H^{n+1})^{>H^{n}})),\ f\neq 0.$$
\end{Prop}
\begin{proof}
By induction on $n$. If $n=0$, then we may apply
\ref {characterize local growth} for $\{ 1\}\subseteq H$, since
$H\neq \{ 1\} $ satisfies the growth axiom.
Notice that $H^1=H^{\# ,\{ 1\} }$ in the notation of \ref {characterize local growth}
and we get the assertion for $n=0$.

For the induction step we assume the assertion for $n$ and
apply \ref {characterize local growth} to the subgroups
$H^n\subsetneq H^{n+1}$ of $G^{\#n+1}$.
Since $H^{n+2}=(H^{n+1})^{\# ,H^n}$ in the notation of
\ref {characterize local growth},
we get the assertion for $n+1$.
\end{proof}


\begin{Lem}\label{growth axiom stable under morphism}
If $U\subseteq H$ are subgroups of $G$ satisfying condition
$(\dagger )$ of Lemma \ref {characterize local growth} and $\psi
:({\bf k}((G)),l)\lra ({\bf k}((G)),l)$ is a morphism, then also
$\psi (U), \psi (H)$ satisfy condition $(\dagger )$.
\end{Lem}
\begin{proof}
Take $h_1\in \psi (H)$ and let $f_1\in {\bf k}((\psi (H)^{>\psi
(U)}))$, $f_1>0$. Then there are $h\in H$ and $f\in {\bf
k}((H^{>U}))$, $f>0$ with $h_1=\psi (h)$ and $f_1= \psi (f)$. By
assumption we have $l(h)=\Log (h)<f$. Therefore $\Log (h_1)=
l(\psi (h))= \psi (l(h))\leq  \psi (f)=f_1$ as desired.
\end{proof}

\section {LE-series constructions.} \label{section:LE}
\sn Let $({\bf k},\log)$ be a logarithmic field. Let $({\bf
k}((G)),l)$ be a field with prelogarithmic section, and let $\psi
$ be a morphism  from $({\bf k}((G )),l)$ to $({\bf k}((G )),l)$,
hence $\psi $ is an order preserving embedding $G\lra G$ such that
$ \psi \circ l=l\circ \psi $. We fix a subgroup $H\neq \{ 1\} $ of
$G$. The goal of this section is to generalize the construction of
the LE series field given in \cite [p. 67-72]{DMM2}. More
precisely, we construct a subfield of ${\bf k}((G))^\EL $ which
contains ${\bf k}((H))$ and which is closed under $\Exp $, $\Log $
and $\psi^{EL} $, assuming natural conditions on $H$ and $\psi
^{EL}$. We shall call this field $LE=LE(H,\psi )$ and suppress the
dependency on $H$ and $\psi $ whenever the data $H$ and $\psi $
are clear from the context (cf. Theorem \ref{LEH}).

\mn
If
$m\in \N _0$ we  write $\psi ^{(m)}$ for $\overbrace {\psi \circ
\cdots\circ \psi }^{m\ times}$. Observe that $\psi ^{(m)} $ is
again a morphism from $({\bf k}((G )),l)$ to $({\bf k}((G )),l)$.
\begin{Def} We define an increasing sequence of subgroups of $G^\EL\>;$
$$G^{0}_m\subseteq G^{1}_m\subseteq \cdots G^{n}_m\subseteq
G^{n+1}_m \cdots
$$
with $G^{n}_m\subseteq G^{\# n}$ as follows:\\
Let $G^{-1}_m:=\{ 1\}$, $G^{0}_m:=\psi ^{(m)}(H)$
which is an ordered multiplicative
subgroup of $G$.
We define by induction on $n$:
$$G^{n+1}_m=G^n_m\cdot \Exp {\bf k}(([G^n_m]^{>G^{n-1}_m}))\subseteq G^{\# n+1}.$$
\sn Set ${\bf k}^n_m:={\bf k}((G^n_m))\>$  (so  ${\bf k}^n_m
\subseteq {\bf k}^{n+1}_m\>$) and $L_m:=\bigcup _n{\bf k}^n_m \>.$
\end{Def}

\sn We now study closure properties of the fields $L_m$ under the
various maps.
\begin{Prop}\label{basic LE}
We have:\sn
 \ (1)\ $\Exp ({\bf k}_m^n)\subseteq {\bf k}_m^n\cdot
G_m^{n+1}\subseteq {\bf k}_m^{n+1}\>$ (so $L_m$ is closed under
$\Exp$).
\sn \ (2) \  $\psi ^\EL (G_m^n)=G_{m+1}^n$.
\end{Prop}
\begin{proof}
(1). We have ${\bf k}_m^n=
{\bf k}(([G_m^n]^{<G_m^{n-1}}))\oplus {\bf k}((G_m^{n-1}))\oplus
{\bf k}(([G_m^n]^{>G_m^{n-1}}))$. Since
$\Exp {\bf k}(([G_m^n]^{<G_m^{n-1}}))\subseteq
{\bf k}((G_m^n))$ we obtain
$$\Exp {\bf k}_m^n\subseteq {\bf
k}_m^n\cdot \Exp {\bf k}((G_m^{n-1}))\cdot
\Exp {\bf k}(([G_m^n]^{>G_m^{n-1}})).$$
By induction hypothesis we know $\Exp {\bf
k}((G_m^{n-1}))\subseteq {\bf k}_m^n$  and by definition of
$G_m^{n+1}$ we have $\Exp {\bf k}(([G_m^n]^{>G_m^{n-1}}))\subseteq
G_m^{n+1}$, so (1) follows. (2). If $n=0$, then $\psi
(G_m^0)=G_{m+1}^0$ by definition. By Proposition \ref {properties
of psi}, $\psi ^\EL $ is a ${\bf k}$-homomorphism which respects
$\Exp $ and arbitrary sums; by induction on $n$, we get $\psi ^\EL
(G_m^n)=G_{m+1}^n$ from the definition of $G_m^n$.
\end{proof}

\begin{Rem}
Notice that it does not follow from the definitions that
$G^{n+1}_m$ is the antilexicographic product of $G^n_m$ and $\Exp
{\bf k}(([G^n_m]^{>G^{n-1}_m}))$.
Similarly, it
does not follow from the definitions that $L_m \subseteq L_{m+1}$.
However, this will be the case under additional assumptions on
$l$, $\psi$ and $H$ which we will introduce step by step in the next statements.
These results will be needed in section \ref{LEinEL}.
\end{Rem}

\begin{Prop}\label{growth axiom in LE}
If $H$ satisfies the growth axiom
(cf. Definition \ref {defn:growth_axiom}) then  for every $m\in \N _0$, $G^0_m$ satisfies the growth axiom
and $G^{n+1}_m$ is the anti-lexicographical product of
$G^n_m$ and
 $\Exp {\bf k}(([G^n_m]^{>G^{n-1}_m}))$.
 Moreover we have
 $$\Exp \kk (( (G^{n}_m)^{>1} ))\subseteq G_m^{n+1}.$$
\end{Prop}
\begin{proof}
$G^0_m$ satisfies the growth axiom by \ref {growth axiom stable under morphism}
applied to $H, \{1\} $ and the morphism $\psi ^{(m)}$.
Now \ref {growth axiom for H} applied to $G^0_m$ shows that
$G^{n+1}_m$ is the anti-lexicographical product of
$G^n_m$ and $\Exp {\bf k}(([G^n_m]^{>G^{n-1}_m}))$.
In particular $G^{n}_m$ is a convex subgroup of $G^{n+1}_m$.
To see the moreover part, take $f\in \kk (( (G^{n}_m)^{>1} ))$.
Since $G^{n}_m$ is a convex subgroup of $G^{n+1}_m$
we can write $f=f_0+...+f_n$, where $\supp f_i\subseteq
(G^{i}_m)^{>G^{i-1}_m}$.
By definition of $G_m^{i+1}$ we have $\Exp f_i\in G_m^{i+1}\subseteq G_m^{n+1}$.
Thus $\Exp f=\prod _{i=0}^n\Exp f_i\in G_m^{n+1}$.
\end{proof}

\begin{Prop}\label{crucial LE}
Assume $H$ satisfies the growth axiom,
 $\psi ^\EL $ is surjective and
\begin{equation} \label{crucialassumption}
 (\psi ^\EL )
^{-1}(G_0^0)\subseteq G_0^1 \>\> (i.e.\>\>(\>\psi ^\EL )
^{-1}(H)\subseteq H\cdot \Exp {\bf k}((H^{>1}))\>) .
\end{equation} Then for all $n,m$ we have:
\sn \ (1) \ $(\psi ^\EL
)^{-1}(G_m^n) \subseteq G_m^{n+1}$.
\sn \ (2) \ ${\bf
k}_m^n\subseteq {\bf k}_{m+1}^{n+1}$.

\medskip
\noindent
If in addition
\begin {equation}\label{crucialassumption2}
l(H)\subseteq \kk((\psi (H))),
\end{equation}
then
\sn \ (3)
$\Log G_m^n\subseteq \mathbf{k}_{m+1}^n$ and
\sn \ (4)
$\Log ({\bf k}_m^n)^{>0}\subseteq {\bf k}_{m+1}^{n+1}$
for all $n,m$.
\end{Prop}
\begin{proof}
(1). We fix $m$ and  show (1) by induction on $n$. Firstly, since
$G_m^0=\psi^{(m)}(H)$, we have $(\psi ^\EL
)^{-1}(G_m^0)=\psi^{(m)}(\psi ^{-1}(H))\subseteq
\psi^{(m)}(G_0^1)=G_{m}^1$, by assumption (\ref {crucialassumption}). Assume we know (1)
for $n$. Then
$$(\psi ^\EL )^{-1}(G_m^{n+1})=(\psi ^\EL
)^{-1}(G_m^n)\cdot (\psi ^\EL )^{-1}(\Exp {\bf
k}(([G^n_m]^{>G^{n-1}_m})))$$  $\subseteq G_m^{n+1}\cdot \Exp
{\bf k}(((\psi ^\EL )^{-1}([G^n_m]^{>G^{n-1}_m})))$ by induction,
so by definition of $G_m^{n+2}$ it remains to show that
$\Exp \kk (( (\psi ^\EL )^{-1}([G^n_m]^{>G^{n-1}_m})))\subseteq G_m^{n+2}.$
By induction we know
$(\psi ^\EL )^{-1}([G^n_m]^{>G^{n-1}_m})\subseteq
(G^{n+1}_m)^{>1}$ and by \ref {growth axiom in LE}
we know $\Exp \kk (( (G^{n+1}_m)^{>1} ))\subseteq G_m^{n+2}$.
Hence the claim follows.

\smallskip
\noindent (2). By  Proposition \ref {basic LE}(2) we know
$G_m^n=(\psi ^\EL )^{-1}(G_{m+1}^n)$. Hence ${\bf k}_m^n={\bf
k}((G_m^n))={\bf k}(((\psi ^\EL )^{-1}(G_{m+1}^n)))$ and the
latter field is contained in ${\bf k}_{m+1}^{n+1}={\bf
k}((G_{m+1}^{n+1}))$ by (1).

\smallskip
\noindent
Now assume that in addition $l(H)\subseteq \mathbf{k}_1^1$.

\smallskip
\noindent (3).
We have $\Log G_m^0=\Log \psi ^{(m)}(H)=
\psi ^{(m)}(\Log H)\subseteq \psi ^{(m)}(\mathbf{k}_1^0)$ by assumption.
Since $\psi ^{(m)}(\mathbf{k}_1^0)=\mathbf{k}_{m+1}^0$ we get
$\Log G_m^0\subseteq \mathbf{k}_{m+1}^0$.
By definition of $G_m^{n+1}$ we have
$\Log (G_m^{n+1})\subseteq \Log G_m^n+{\bf k}_m^n$.
By (2) we know $\kk _m^n\subseteq {\bf k}_{m+1}^{n+1}$
and
by induction on $n$ we have
$\Log G_m^n \subseteq {\bf k}_{m+1}^{n}$.
Since ${\bf k}_{m+1}^{n}\subseteq {\bf k}_{m+1}^{n+1}$
we obtain (3).

\smallskip
\noindent (4).
Since $({\bf k}_m^n)^{>0}=G_m^n\cdot {\bf
k}^{>0}\cdot (1+{\bf k}(((G_m^n)^{<1})))$ we get $\Log ({\bf
k}_m^n)^{>0}\subseteq \Log (G_m^n)+{\bf k}((G_m^n))$.
Now (2) says ${\bf k}((G_m^n))
\subseteq {\bf k}_{m+1}^{n+1}$ and (3) says
$\Log (G_m^n)\subseteq {\bf k}_{m+1}^{n}$.
Since ${\bf k}_{m+1}^{n}\subseteq {\bf k}_{m+1}^{n+1}$ we obtain (4).
\end{proof}

Hence under the assumptions of Proposition \ref{crucial LE}, we have
$L_0\subseteq L_1\subseteq \cdots$ and we define $\LE:=\LE (H\>;
\psi):=\bigcup L_m$.
Here the situation in a
diagram where a line indicates containment:


$$\vcenter
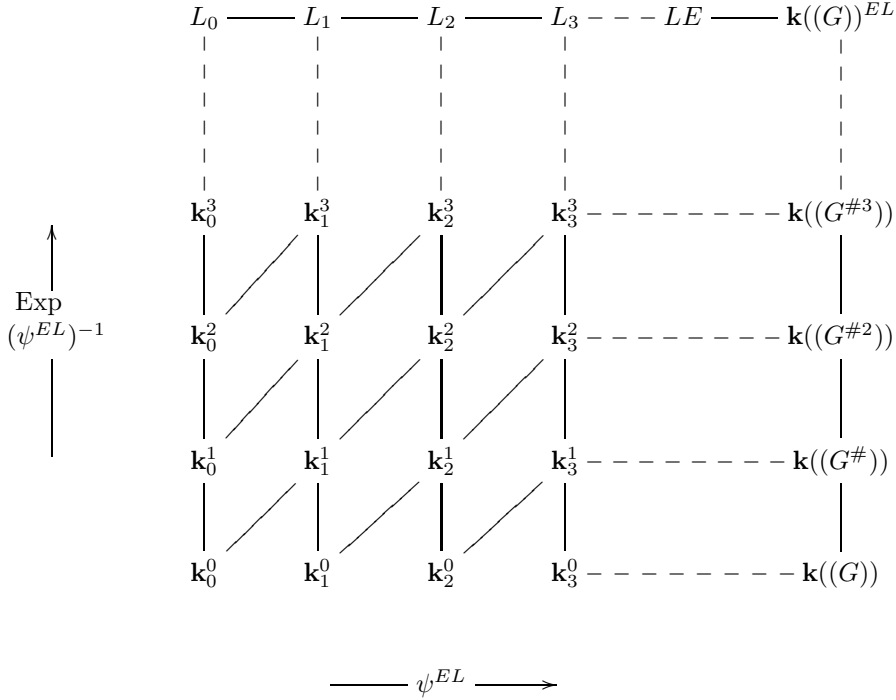
{\xymatrix@=25pt{
&L_0\ar@{-}[r]&        L_1\ar@{-}[r]&         L_2\ar@{-}[r]&        L_3\ar@{--}[r]&LE\ar@{-}[r]        &\mathbf{k}((G))^\EL \\
&&&&&&\\
&\mathbf{k}^3_0\ar@{--}[uu]&\mathbf{k}^3_1\ar@{--}[uu]&\mathbf{k}^3_2\ar@{--}[uu]&\mathbf{k}^3_3\ar@{--}[uu]\ar@{--}[rr]&&\mathbf{k}((G^{\# 3}))\ar@{--}[uu]\\
\vbox{\hbox{\ \ $\Exp $} \hbox{\ $(\psi ^\EL )^{-1}$}}&\mathbf{k}^2_0\ar@{-}[u]\ar@{-}[ur]&\mathbf{k}^2_1\ar@{-}[u]\ar@{-}[ur]&\mathbf{k}^2_2\ar@{-}[u]\ar@{-}[ur]&\mathbf{k}^2_3\ar@{-}[u]\ar@{--}[rr]&&\mathbf{k}((G^{\# 2}))\ar@{-}[u]\\
{\ }\ar '[u][uu]&\mathbf{k}^1_0\ar@{-}[u]\ar@{-}[ur]&\mathbf{k}^1_1\ar@{-}[u]\ar@{-}[ur]&\mathbf{k}^1_2\ar@{-}[u]\ar@{-}[ur]&\mathbf{k}^1_3\ar@{-}[u]\ar@{--}[rr]&&\mathbf{k}((G^\# ))\ar@{-}[u]\\
&\mathbf{k}^0_0\ar@{-}[u]\ar@{-}[ur]&\mathbf{k}^0_1\ar@{-}[u]\ar@{-}[ur]&\mathbf{k}^0_2\ar@{-}[u]\ar@{-}[ur]&\mathbf{k}^0_3\ar@{-}[u]\ar@{--}[rr]&&\mathbf{k}((G))\ar@{-}[u]\\
&&{\ } \ar '[r][rr]&\psi ^\EL &&&\\
}}$$

So far we know that
$\LE$ is closed under $\Exp $, $\Log $
and $\psi^{EL} $.

\begin{Thm}\label{LEH} Under the assumptions of Proposition \ref{crucial LE}, we have
$L_0\subseteq L_1\subseteq \cdots$ and we define $\LE:=\LE (H\>;
\psi):=\bigcup L_m$ - thus $\LE$ is closed under $\Exp $, $\Log$
and $\psi^{EL} $.
\end{Thm}

%
\noindent We will now examine a sufficient condition on $l$, $\psi
$ and $H$ for assumption (\ref{crucialassumption}) of Proposition
\ref {crucial LE} to be fulfilled:
$$\psi \vert _{l^{-1}(G)\cap H}=l\vert _{l^{-1}(G)\cap H}\leqno ({\rm Comp}_H).$$
Notice that $l^{-1}(G)=l^{-1}(G^{>1})$ as $l$ has values in $\mathbf{k}((G^{>1}))$.

\begin{Lem}\label{compute psiE first level}
Let $l$, $\psi$, $H$ be such that $({\rm Comp}_H)$ holds. Then for
every $h\in l^{-1}(G)\cap H$, and for every $k \in {\bf k}$ we
have: \sn (1) $\psi ^\EL (\Exp (k\cdot h))=\Exp (k \cdot \Log
(h))$. \sn  (2) If $\psi ^\EL $ is surjective, then $(\psi ^\EL
)^{-1} (\Exp (k \cdot \Log (h)))\in G_0^1 \>.$
\end{Lem}
\begin{proof}
(1). Take $h\in l^{-1}(G)\cap H$. Then $\psi ^\EL  (\Exp (k \cdot
h))= \Exp (\psi ^\EL (k \cdot h))=\Exp (k \cdot \psi ^\EL (h))=
\Exp (k \cdot \psi (h))=\Exp (k \cdot \Log (h))$. (2). By (1),
$(\psi ^\EL )^{-1}(\Exp (k \cdot \Log (h)))=\Exp (k \cdot h)$.
Since $h\in l^{-1}(G)\cap H$ we have $h>1$. Thus $\Exp (k \cdot
h)\in G_0^1$ by definition.
\end{proof}

\begin{Cor}\label{scholium}
Let $l$, $\psi$, ${H_0}$ be such that $({\rm Comp}_{H_0})$ holds, and
assume that $\psi ^\EL $ is surjective. Let $H$ be
a subgroup of $G$ containing ${H_0}$ and contained in the group
generated by all the $\Exp (k \cdot \Log (h))$, $k \in {\bf k}$,
$h\in l^{-1}(G)\cap {H_0}$. Then assumption (\ref {crucialassumption}) of Proposition \ref
{crucial LE} is satisfied for $H$.
\end{Cor}
\begin{proof}
We must show that $(\psi ^\EL )^{-1}$ maps $H$ into
$G_0^1=H\cdot\Exp ({\bf k}((H^{>1})))$.
By assumption on $H$ it suffices to show that
$(\psi ^\EL )^{-1}$ maps each
element $\Exp (k \cdot \Log (h))$, $k \in {\bf k}$,
$h\in l^{-1}(G)\cap {H_0}$ into $G_0^1$.
As $({\rm Comp}_{H_0})$ is satisfied we may apply
\ref {compute psiE first level}(2) for $H_0$, which precisely says this.
\end{proof}


\section {Finding the LE-series field in the exponential field generated
by logarithmic words.}\label{LEinEL}\sn
 Let $G$ be the
multiplicative group of logarithmic words: $$G:=\>\{x^{r _0}\cdot
(\log x)^{r _1}\cdot ...\cdot (\log _nx)^{r _n}, n\in \N _0, r
_i\in k \}\>.$$ Define
$$\psi (x^{r _0}\cdot (\log x)^{r _1}\cdot ...\cdot (\log _nx)^{r _n}):=
(\log x)^{r _0}\cdot (\log _2x)^{r _1}\cdot ...\cdot (\log
_{n+1}x)^{r _n}.$$ So $\psi :G\lra G$ is an order preserving group
embedding. Let $l\>:\> G \lra {\bf k}((G ^{>1}))$ be defined by
$$l(x^{r _0}\cdot (\log x)^{r _1}\cdot ...\cdot (\log _nx)^{r
_n}):= r _0\cdot \log x \>+\> \cdots \>+\> r _n\cdot \log _{n+1}x
\>.$$ Then \sn \
(1) \ $l$ is a prelogarithmic section of ${\bf
k}((G ))$ and $\psi $ is a morphism from $({\bf k}((G)),l)$ to
$({\bf k}((G )),l)$. This is obvious.\sn \
(2)\ $\psi ^\EL $ is
surjective, by Proposition \ref {properties of psi}(3) (applied to
$n=1$).\sn \
(3) \ Clearly $\psi$ is a contracting morphism. Hence
by Proposition \ref {psi contracting}, $\psi ^{\EL}$ is a contraction.
In particular $(\psi ^{EL})^{-1}(g)>g$ for all $g\in G^{\EL}$,
$g>1$. \sn \
(4) \ $H:=\{ x^k \st k \in {\bf k} \} $ is a subgroup
of $G$ and $H_0:=\{ x^z\st z\in \Z \} $ is a subgroup of $H$,
generated by $l^{-1}(G)\cap H_0=l^{-1}(G)\cap H=\{ x\}$.
\sn \
(5) \ $H$
satisfies the growth axiom
  (cf. Definition \ref {defn:growth_axiom}), since
$l(x^k )=k \log x<\vert f\vert $ for all $f\in {\bf k}((H^{>1}))$,
$k\in {\bf k}$. \sn \
(6) \ $\psi $ satisfies ${\rm Comp}_H$,
since $l^{-1}(G)\cap H=\{x \}$ and $\psi (x)=l(x)$.
\sn \
(7) \ $l$ and $\psi $ satisfy condition (\ref{crucialassumption2}) (cf. \ref{crucial LE}) as $l(x^k)=k\log x=k\psi (x)\in \mathbf{k}((\psi (H)))$

\mn Conditions (2),(3) and (6) show that
all assumptions of \ref {scholium} are satisfied.
By this property (\ref{crucialassumption}) of \ref{crucial LE} and by (7) also property
(\ref{crucialassumption2}) of \ref{crucial LE} is satisfied.
Hence the
LE-series construction from section \ref {section:LE} is
applicable and we show that it gives the field ${\bf
k}((x^{-1}))^{LE}$ of LE-series constructed by \cite[p. 67-72]
{DMM2}, where the composition with $\Exp (x)$ is
$$\Phi :=(\psi ^\EL )^{-1}.$$
By induction on $n$ Proposition \ref {growth axiom in LE}
identifies the group $G_0^n$ with the group $G_n$ from the
construction in \cite {DMM2}. Notice that by induction on $n$ we
may use our exponential function from ${\bf k}((G))^{EL}$ as the
abstract isomorphism taken in \cite {DMM2} to define their
exponential in step $n$. We see that $L_0$ is the field ${\bf
k}((x^{-1}))^E$ from \cite {DMM2}. Moreover we have shown that
$\Phi $ satisfies the following properties:
\begin{enumerate}
\item $\Phi $ is an order preserving ${\bf k}$-isomorphism of
exponential ordered fields. \item $\Phi (x^k )=\Exp(k \cdot x)$
$(k \in {\bf k})$ \item  $\Phi $ respects arbitrary sums. \item
$L_m\subseteq L_{m+1}$ \item $\Phi \vert _{L_{m+1}}:L_{m+1}\lra
L_m$ is an isomorphism.
\end{enumerate}
\medskip
By induction on $m$, these properties identify the exponential
fields $L_m$ with the fields $L_m$ from \cite {DMM2} and the
composition  with $\Exp (x)$ (defined on $L_m$) as $(\psi ^\EL
)^{-1}$ (restricted to $L_m$). This shows that the exponential
field ${\bf k}((x^{-1}))^{\LE}$ is  our exponential field $LE$ and
the composition with $\Exp (x)$ is $(\psi ^\EL )^{-1}$.

\section { The ${\bf L} $-operation.} \sn \label{notELinLE}
Let $\Gamma $ be a totally ordered set, $\sigma :\Gamma \lra
\Gamma $ a decreasing embedding, $({\bf k}, \log)$ an ordered
field with a surjective logarithm. \sn In this section, we show
that the exponential-logarithmic field ${\bf k} ((\Gamma^{\bf k}))^{\sigma \EL}\>\>$ (cf.\ Examples \ref{sigmamultiplicativehahn}
and \ref{sigmamultiplicativehahn2} and \cite{KS}) does not embed
in the logarithmic-exponential field ${\bf k}((x^{-1}))^{\LE}\>\>$
(cf.\ \cite {DMM2}). We need the following notions:
\begin{Def}
Let $(K_1, v_1)$ and $(K_2, v_2)$ be valued fields with value
groups $G_1$ and $G_2$ respectively. An embedding $\phi : K_1 \lra
K_2$ is said to be {\it valuation preserving} if the induced map
$v_1(a) \mapsto v_2(\phi(a))$ is a well defined order preserving
embedding of $G_1$ into $G_2$.
\end{Def}
\begin{Def}
As in \cite[6.27]{DMM2} we introduce the following notion. Fix
$\alpha\in {\bf k} ((\Gamma^{\bf k}))^{\sigma \EL}\>\>$ such that $\alpha
>k$ for all $k\in {\bf k}$. Let $z \in \Z$. Let $s\in {\bf k} ((\Gamma^{\bf k}))^{\sigma
\EL}$ such that $s
>k$ for all $k\in {\bf k}$. We say that $s$ has {\it level} $z$
with respect to $\alpha$ if there exists $N\in N_0$ such that
$v(\Log _{N+z}(s))=v(\Log_N(\alpha))$.
\end{Def}
\def\phi{\varphi}
\begin{Def}
As in \cite{DMM2}, the exponential on ${\bf k} ((x^{-1}))^{\LE}$
will be denoted by $E$ and the logarithm by $L$. An embedding
$\phi :{\bf k} ((\Gamma^{\bf k}))^{\sigma \EL}\lra {\bf k}
((x^{-1}))^{\LE}$ is said to be {\it exponential preserving} if it
commutes with the exponential maps defined on its domain and
co-domain, i.\ e.\ if $E(\phi( \alpha)) = \phi (\Exp(\alpha))$ for
all $\alpha \in {\bf k} ((\Gamma^{\bf k}))^{\sigma \EL}$. Note that in
this case, $\phi$ will also be {\it logarithm preserving}, that
is, $L(\phi( \alpha)) = \phi (\Log(\alpha))$ for all $\alpha \in
{\bf k} ((\Gamma^{\bf k}))^{\sigma \EL}$, $\alpha >0$.
\end{Def}
\begin{Thm}
\label{base cover production}
 Then there is no
exponential $\bf k$--embedding of ordered fields from $\>{\bf k}
((\Gamma^{\bf k}))^{\sigma \EL}$ into ${\bf k} ((x^{-1}))^{\LE}\>.$
\end{Thm}
\begin{proof}
Fix an element $X:=x_{\gamma}$ for some $\gamma \in \Gamma $ and
consider
$${\bf L} :=\sum _{n\geq 0}\Log _{\sigma ,n} X.$$ Note that by definition,
$\Log _{\sigma ,n} X = \sigma_n (X)$ and the support of ${\bf L}$
is anti wellordered, so ${\bf L} $ is a well defined element of
${\bf k} ((\Gamma^{\bf k}))^{\sigma \EL}$. \sn Now suppose that there is
an exponential $\bf k $--embedding $\phi :\>{\bf k} ((\Gamma^{\bf k}))^{\sigma \EL}\rightarrow {\bf k} ((x^{-1}))^{\LE}.$ Note that
$\phi$ is valuation preserving since it is a ${\bf k}$--embedding.
Let $Y:=\phi (X)\in {\bf k} ((x^{-1}))^{\LE}$. Then $Y>{\bf k}$
and by \cite[Prop. 6.3]{DMM2} there is a compositional inverse
$Y^{-1}$ of $Y$ in ${\bf k} ((x^{-1}))^{\LE}$. By \cite[Th.
6.2]{DMM2} that means there is an $E$--preserving ${\bf
k}$--embedding $\psi :{\bf k}((x^{-1}))^{\LE}\lra {\bf k}
((x^{-1}))^{\LE}$ of ordered fields with $\psi (Y)=x$. So, by
replacing $\phi $ with $\psi \circ \phi $ we may assume that $\phi
(X)=x$.  By \cite[Lemma 6.21]{DMM2}, there is some $i\in \N $ such
that $\Phi ^i(\phi ({\bf L} ))\in L_0$. As $\Phi ^i:{\bf k}
((x^{-1}))^{\LE}\lra {\bf k} ((x^{-1}))^{\LE}$ is an $\exp
$--preserving {\bf k}--embedding  of ordered fields, we may
replace $\phi $ by $\Phi ^i\circ \phi $ and assume that $\phi
(X)=E _i(x)$ and $\phi ({\bf L} )\in L_0$. Let ${\bf T}:={\bf L} -
\sum_{n=0}^{i-1}\Log _{\sigma ,n}X$ and ${\bf S}:=\Log _{\sigma
,i}X$. Then $\phi ({\bf S})=x$ and $\phi ({\bf T})\in L_0$. Now
$\phi$ is level preserving since it is both log-preserving and
valuation preserving. Consider ${\bf T} -{\bf S}>{\bf k} $. We
compute that $\mbox{ level }({\bf T} - {\bf S})< \mbox{ level }
({\bf S})$ (w.r.t. ${\bf S}$). Hence $\mbox{ level }(\phi ({\bf
T})-x)<\mbox{ level }(x)=0$ (w.r.t. $\phi ({\bf S})=x$). So $\phi
({\bf T})-x\in L_0$ is an element $>{\bf k} $ having a level $<0$.
This contradicts \cite[Prop. 6.28 p. 100]{DMM2}).\end{proof} \bn


\def\lit #1#2#3#4{%
\bibitem [#1]{#1} #2 \textit{#3;} #4
}%

\end{document}